\documentclass[11pt]{amsart}

\usepackage{amsmath}
\usepackage{amssymb}
\usepackage{mathrsfs}

\usepackage{epsf}

\newtheorem{lemma}{Lemma}[section]
\newtheorem{theorem}[lemma]{Theorem}
\newtheorem{fact}[lemma]{Fact}
\newtheorem{definition}[lemma]{Definition}
\newtheorem{corollary}[lemma]{Corollary}
\newtheorem{proposition}[lemma]{Proposition}
\newtheorem{conjecture}[lemma]{Conjecture}

\newcommand{\opt}{\textrm{opt}}
\newcommand{\grd}{\textrm{grd}}
\begin{document}

\author[Anna Niewiarowska, Michal Adamaszek]{Anna Niewiarowska,
Michal Adamaszek$^*$}
\thanks{$^*$ University of Warsaw, Dept. of Mathematics, Informatics and Mechanics,
ul. Banacha 2, 02-097 Warsaw, Poland, email: \texttt{\{annan,aszek\}@mimuw.edu.pl}}
\title{Combinatorics of the change-making problem}

\begin{abstract}
We investigate the structure of the currencies (systems of coins) for which the
greedy change-making algorithm always finds an optimal solution (that is, a one
with minimum number of coins). We present a series of necessary conditions that
must be satisfied by the values of coins in such systems. We also uncover
some relations between such currencies and their sub-currencies.
\end{abstract}


\maketitle

\section{Introduction}
In the change-making problem we are given a set of coins and we wish to
determine, for a given amount $c$, what is the minimal number of coins needed
to pay $c$. For instance, given the coins $1,2,5,10,20,50$, the minimal
representation of $c=19$ requires $4$ coins ($10+5+2+2$).

This problem is a special case of the general knapsack problem with all coins
of unit weights. In some cases the solution may be found by a greedy strategy
that uses as many of the largest coin as possible, then as many of the next one
as possible and so on. This greedy solution is optimal for the set of
coins given above, but fails to be optimal in general. For instance, if we have
coins $1,5,9,16$, then the amount $18$ will be paid greedily as $16+1+1$, while
the optimal solution ($9+9$) requires just two coins. In this paper we shall
concentrate on the combinatorial properties of those sets of coins for which
the greedy solution is always optimal.

The sequence $A=(a_0,a_1,\ldots,a_k)$, where $1=a_0<a_1<\ldots<a_k$ will be
called a {\it currency} or {\it coinage system}. We always assume $a_0=1$ is
the smallest coin to avoid problems with non-representability of certain
amounts. For any amount $c$ by $\opt_A(c)$ and $\grd_A(c)$ we denote,
respectively, the minimal number of coins needed to pay $c$ and the number of
coins used when paying $c$ greedily (for example, if $A=(1,5,9,16)$ then
$\opt_A(18)=2$ and $\grd_A(18)=3$). The currency $A$ will be called {\it
orderly}\footnote{Various authors have used the terms orderly
\cite{Jones,Maurer}, canonical \cite{KoZa,Pear}, standard \cite{TienHu} or
greedy \cite{CCS}.}
if for all amounts $c>0$ we have $\opt_A(c)=\grd_A(c)$.

If a coinage system $A$ is not orderly then any amount $c$ for which
$\opt_A(c)<\grd_A(c)$ will be called a {\it counterexample}.

Let us briefly summarize related work. Magazine, Nemhauser and Trotter
\cite{MNT} gave a necessary and sufficient condition to decide whether
$A=(1,a_1,\ldots,a_{k+1})$ is orderly provided we know in advance that
$A'=(1,a_1,\ldots,a_k)$ is orderly, the so-called one-point theorem (see
section \ref{section2} of the present paper). Kozen and Zaks \cite{KoZa} proved,
among other things, that the smallest counterexample (if exists), does not
exceed the sum of the two highest coins. They also asked if there is a
polynomial-time algorithm that tests if a coinage system is orderly. Such an
algorithm was presented by Pearson \cite{Pear}. It produces a set of $O(k^2)$
``candidates for counterexamples'', which is guaranteed to contain the smallest
counterexample if one exists. The rest of the algorithm is just testing these
potential candidates, and the overall complexity is $O(k^3)$. A similar set of
possible counterexamples (perhaps not containing the smallest one), but of size
$O(k^3)$, was given by Tien and Hu in \cite{TienHu} (see formula (4.20) and
Theorem 4.1 of that paper). It leads to an $O(k^4)$ algorithm. The authors of
\cite{MNT} and \cite{TienHu} were concentrating mainly on the error analysis
between the greedy and optimal solutions. Apparently Jones \cite{Jones} was the
only one who attempted to give a neat combinatorial condition characterizing
orderly currencies, but his theorem suffered from a major error, soon pointed
out by Maurer \cite{Maurer}. Our paper has been paralleled by an independent
work of Cowen, Cowen and Steinberg \cite{CCS} about currencies whose all
prefixes are orderly and about non-orderly currencies which cannot be
``fixed'' by appending extra coins. The results contained in section
\ref{section4} and a special case ($l=2$) of Theorem \ref{theoremprefixorderly}
of this paper have also been proved in \cite{CCS}.

The aim of this paper is to study some orderly coinage systems from a
combinatorial viewpoint, motivated by the need to have some
nice characterization. One of the motivations was the observation that 
\emph{if $A=(1,a_1,a_2,\ldots,a_k)$ is an orderly currency, then the currency
$(1,a_1,a_2)$ is also orderly}, that will be generalized and proved in Theorem
\ref{theoremprefixorderly}. Going further, one may start with an orderly
currency, take out some of its coins and ask if the remaining coins again form
an orderly currency. The precise answer to this question, given in section
\ref{section7}, will be a consequence of the results of sections \ref{section3}
and \ref{section5}, where we prove some properties of the distances $a_j-a_i$
between the coins of an orderly currency. In section \ref{section4} these
results will be used to give a complete description of orderly currencies with
less than 6 coins. In section \ref{section6} we study the behaviour of the
currencies obtained as prefixes of an orderly currency. Some closing remarks and
open problems are included in section \ref{section8}.

\section{Preliminary results}
\label{section2}

If $A=(1,a_1,\ldots,a_k)$ is a currency it will often be convenient to set
$a_{k+1}=\infty$. This will be especially useful whenever we want to choose,
say, the first interval $[a_m,a_{m+1}]$ of length at least $d$ for some $d$.
The reader will see that in all applications the infinite interval will be as
useful as proper intervals.

There are three standard arguments that will be used repeatedly throughout this
paper, so we quote them now to avoid excessive repetitions in the future.
All the time we assume $A=(1,a_1,\ldots,a_k)$ is orderly.

First, suppose we have $a_m$, $a_i$ and $a_l$, such that $a_l<a_m+a_i<a_{l+1}$.
Then $a_m+a_i$ has a representation that uses 2 coins. Since $a_m+a_i$ is
strictly between $a_l$ and $a_{l+1}$ its greedy representation must start
with $a_l$, followed (since $A$ is orderly) by just one other coin $a_r$.
It follows that there exists $r$ such that $a_m+a_i=a_l+a_r$.
 
The second argument is a slight modification of the first one; namely, if
$a_l<a_m+a_i$ and the number $a_m+a_i-a_l$ is not one of the coins, then
$a_{l+1}\leq a_m+a_{i}$.

The third argument is a bit more complicated. Suppose that for
some $j>i\geq 1$ we have $a_j-a_i=d$. Let us choose the largest $m$ for which
$a_m-a_{m-1}<a_i$ (such $m$'s exist, for instance $a_i-a_{i-1}<a_i$). Then
$a_{m+1}-a_m\geq a_i$ (here it is possible that $a_{m+1}=\infty$), so we
have $a_m<a_{m-1}+a_i<a_{m+1}$. If we also have $a_m<a_{m-1}+a_j<a_{m+1}$ then,
as before, there exist numbers $r<i$ and $s<j$ such that:

\begin{center}
\begin{tabular}{c}
$a_{m-1}+a_i=a_m+a_r,$\\
$a_{m-1}+a_j=a_m+a_s.$
\end{tabular}
\end{center}
\begin{figure}
\epsfbox{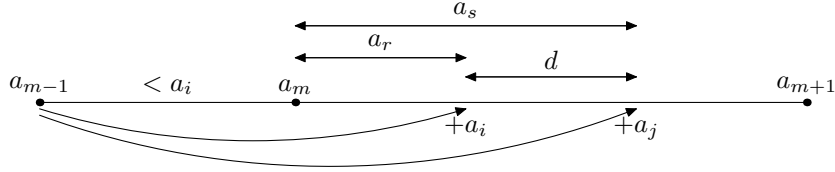}
\caption{Illustration of a standard argument}
\end{figure}
Then $a_s-a_r=a_j-a_i=d$, so we decreased the coins' indicies from $(j,i)$ to
$(s,r)$, keeping the difference $d$ unchanged. Therefore, if additionally
$(j,i)$ was the {\it smallest} pair of indices for which $a_j-a_i=d$, we would
have a contradiction, hence we may assume that in such case $a_{m-1}+a_j\geq
a_{m+1}$.

We shall frequently make use of the following famous result:

\begin{theorem}[{\bf One-point theorem, \cite{MNT,HuLen,CCS}}]
\label{onept}
Suppose $A'=(1,a_1,\ldots,a_k)$ is orderly and $a_{k+1}>a_k$. Let $m=\lceil
a_{k+1}/a_{k}\rceil$. Then $A=(1,a_1,\ldots,a_k,a_{k+1})$ is orderly if and
only if $\opt_A(ma_k)=\grd_A(ma_k)$.
\end{theorem}

{\bf Remark.} According to this theorem, if the shorter currency $A'$ is
orderly, then the optimality of the greedy solution for $A$ needs to be checked
only for the single value $ma_k$. This justifies the name {\it one-point
theorem}. Note, that although in general it is NP-hard to compute $\opt_A(c)$
for arbitrary $A$ and $c$ (see \cite{Shallit} and \cite{KoZa} for a discussion),
the one-point theorem test $\opt_A(ma_k)=\grd_A(ma_k)$ runs in polynomial time,
since it is equivalent to $\grd_{A'}(ma_k-a_{k+1})\leq m-1$. For the sake of
completeness we decided to include a short proof of the one-point theorem.

{\bf Proof.} One of the implications is trivial. Now suppose that
$\opt_A(ma_k)=\grd_A(ma_k)$. We have
$$(m-1)a_k+1\leq a_{k+1}\leq ma_k.$$
For all values $c< a_{k+1}$ all the payments $\grd_A(c)$, $\grd_{A'}(c)$,
$\opt_A(c)$, $\opt_{A'}(c)$ coincide, so $\grd_A(c)=\opt_A(c)$. All other $c$
will be split in two groups: $c\in[a_{k+1}, ma_k)$ and $c\geq ma_k$.

{\bf 1. $a_{k+1}\leq c<ma_k$.} For every such $c$ we have $c<2a_{k+1}$,
therefore any payment of $c$ contains either $0$ or $1$ copies of $a_{k+1}$.
Together with the orderliness of $A'$ this implies
$$\opt_A(c)=\min\{1+\grd_{A'}(c-a_{k+1}), \grd_{A'}(c)\}.$$
At the same time $1+\grd_{A'}(c-a_{k+1})=\grd_A(c)$, so in order to prove
$\opt_A(c)=\grd_A(c)$ it suffices to show the inequality
$$\grd_A(c)\leq \grd_{A'}(c).$$
Observe that
$$\grd_{A'}(c)=(m-1)+\grd_{A'}(c-(m-1)a_k).$$
The function $\grd_{A'}=\opt_{A'}$ satisfies the triangle inequality, so
$$\grd_{A'}(ma_k-a_{k+1})+\grd_{A'}(c-(m-1)a_k)\geq
\grd_{A'}(c-a_{k+1}+a_k)=1+\grd_{A'}(c-a_{k+1}).$$
Finally
$$\grd_{A'}(c)-\grd_{A}(c)=(m-1)+\grd_{A'}(c-(m-1)a_k)-(1+\grd_{A'}(c-a_{k+1}
))\geq$$
$$\geq m-2 + 1 -\grd_{A'}(ma_k-a_{k+1})=m-1-\grd_{A'}(ma_k-a_{k+1}).$$
However
$$1+\grd_{A'}(ma_k-a_{k+1})=\grd_A(ma_k)=\opt_A(ma_k)\leq m,$$
which eventually implies the desired inequality
$$\grd_{A'}(c)-\grd_{A}(c)\geq 0.$$

{\bf 2. $c\geq ma_k$.} Denote by $\mathcal{OPT}(c)$ the set of optimal payments
for $c$:
$$\mathcal{OPT}(c) = \{(x_0,\ldots,x_{k+1}): \sum_{i=0}^{k+1}x_ia_i=c \textrm{
and } \sum_{i=0}^{k+1} x_i \textrm{ is minimal}\}$$

It is sufficient to exhibit a payment $(x_i)\in\mathcal{OPT}(c)$ with
$x_{k+1}>0$. Consider any optimal payment $(x_i)$. We may apply to it the
following two operations:
\begin{itemize}
\item if $x_k\geq m$ then replace $m$ coins $a_k$ with the greedy decomposition
of $ma_k$. This way the number of coins in the payment does not increase (since
$\opt_A(ma_k)=\grd_A(ma_k)$), while the multiplicity of $a_k$ in the payment
decreases.
\item if $\sum_{i=0}^{k-1}x_ia_i\geq a_k$ then instead of the coins needed to
pay $\sum_{i=0}^{k-1}x_ia_i$ insert the greedy decomposition of this amount
with respect to $A'$. This will not increase the overall number of coins (since
$A'$ was orderly), but it will decrease the amount paid with
$1,a_1,\ldots,a_{k-1}$.
\end{itemize}
It is clear that repeating these two steps sufficiently many times we will
finally end up with an optimal payment $(x_i)$ satisfying
$\sum_{i=0}^{k-1}x_ia_i<a_k$ and $x_k<m$. Then
$$\sum_{i=0}^{k}x_ia_i\leq a_k-1+(m-1)a_k=ma_k-1<c$$
hence $x_{k+1}>0$ in this payment. \qed

It is obvious that the one-coin currency $A=(1)$ is orderly, as well as
all the two-coin currencies $A=(1,a_1)$. The reader may now wish to solve
the easy problem of when a three-coin currency $A=(1,a_1,a_2)$ is
orderly. For reasons which will become clear later we shall express the
solution in terms of the following set:

\begin{definition}
\label{defA}
For any $a>0$ we define:
$$\mathcal{A}(a)=\bigcup_{m=1}^{\infty} \bigcup_{l=0}^m \{ma-l\}=$$
$$=\{a-1,a\}\cup\{2a-2,2a-1,2a\}\cup\ldots\cup\{ma-m,\ldots,ma\}\cup\ldots$$
\end{definition}

\begin{proposition}
\label{lemma3orderly}
The currency $A=(1,a_1,a_2)$ is orderly if and only if
$a_2-a_1\in\mathcal{A}(a_1)$.
\end{proposition}

{\bf Proof.} Let $m=\lceil a_2/a_1\rceil$. By the one-point theorem $A$ is
orderly if and only if the greedy algorithm is optimal for $ma_1$, which is
equivalent to
$$\grd_A(ma_1)\leq m$$
or
$$ma_1-a_2\leq m-1$$
which means that $a_2-a_1 = (m-1)a_1 - (ma_1-a_2) \in \mathcal{A}(a_1)$ (more
precisely, $a_2-a_1$ belongs to the $(m-1)$-st summand of $\mathcal{A}(a_1)$).
On the other hand, if $m$ is the least number for which $a_2-a_1$ belongs to the
$(m-1)$-th summand of $\mathcal{A}(a_1)$, then $\lceil a_2/a_1\rceil=m$ and
$a_2-a_1=(m-1)a_1-l$ for some $l\leq m-1$. Then
$$\grd_A(ma_1)=1+(ma_1-a_2) = 1+l\leq m$$
as required.\qed

\section{Investigating differences, part I}
\label{section3}

In this section we begin investigating distances between the coins of an
orderly coinage system, followed by an application of these results. 

\begin{proposition}
\label{lemmanot1}
If $A=(1,a_1,\ldots,a_k)$ is orderly and $a_1\geq 3$, then 
$$a_{i}-a_{i-1}\neq 1$$
for all $i=1,\ldots,k$.
\end{proposition}

{\bf Proof.} Suppose, on the contrary, that $a_j-a_{j-1}=1$ and let $j$ be the
least index with this property. Since $a_1\geq 3$, we have $j\geq 2$.

Let us choose the largest index $m$ for which $a_m-a_{m-1}<a_{j-1}$. Then
$a_{m+1}-a_m\geq a_{j-1}$, and if $a_{m-1}+a_j<a_{m+1}$ then we would have a
contradiction by the third standard argument from section \ref{section2}.
Therefore $a_{m+1}\leq a_{m-1}+a_j$. Since
$a_{m+1}-a_m\geq a_{j-1}$, we have
$$a_j=a_{j-1}+1\leq (a_{m+1}-a_m)+(a_m-a_{m-1})=a_{m+1}-a_{m-1}\leq a_j$$
meaning that $a_m-a_{m-1}=1$ and $a_{m+1}-a_m=a_{j-1}$.

It follows that
$$a_m<a_{m-1}+a_{j-1}<a_{m+1}$$
which means that $a_{m-1}+a_{j-1}-a_m=a_{j-1}-1$ must be one of the coins,
contradicting the minimality of $j$. This ends the proof.\qed

The previous proposition can be sharpened as follows:
 
\begin{proposition}
\label{lemmadiffa1-1} 
If $A=(1,a_1,\ldots,a_k)$ is orderly then
$$a_{i}-a_{i-1} \geq a_1-1$$
for all $i=1,\ldots,k$.
\end{proposition}

{\bf Proof.} This is obviously true if $a_1=2$, so let $a_1\geq 3$. Let $j$ be
the largest index for which $a_j-a_{j-1}\leq a_1-2$. By
Proposition \ref{lemmanot1} we
have $a_j-a_{j-1}\geq 2$. From the maximality of $j$ we have $a_{j+1}-a_j\geq
a_1-1$ (it is possible that $a_{j+1}=\infty$). Now consider the amount 
$$c=a_{j-1}+a_1.$$
\begin{figure}
\epsfbox{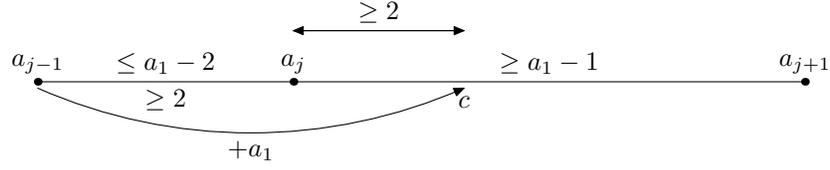}
\caption{Illustration of the proof of Prop. \ref{lemmadiffa1-1}}
\end{figure}
It satisfies $a_j+2\leq c \leq a_j+a_1-2 < a_{j+1}$, hence $\opt_A(c)=2$.
When paid greedily, the amount $c$ is decomposed to $a_j$ and $c-a_j$
copies of the coin $1$, which makes
$$1+(c-a_j)\geq 1+2=3$$
coins altogether. This contradicts the fact $A$ is orderly, thus completing
the proof of this proposition.\qed

Proposition \ref{lemmadiffa1-1} imposes certain restrictions on the possible
differences $a_i-a_{i-1}$. In the next theorem we shall generalize this
restriction, but first let us state without proof some obvious properties of the
sets $\mathcal{A}(a)$ that will be useful in the proof:

\begin{fact} Let $a\geq 2$ be an integer. Then:
\label{factabouta}
\begin{itemize}
\item[(1)] if $x,y\in\mathcal{A}(a)$ then $x+y\in\mathcal{A}(a)$.
\item[(2)] an integer $x\geq 2$ does not belong to $\mathcal{A}(a)$ if and only
if there exists an integer $p\geq 0$ such that
$$pa+1\leq x\leq (p+1)a-(p+2).$$
\item[(3)] if $p_1<p_2<\ldots<p_m$ and $p_m-p_{1}\not\in\mathcal{A}(a)$ then
$p_j-p_{j-1}\not\in\mathcal{A}(a)$ for some $2\leq j\leq m$ (this follows from
(1)).
\item[(4)] if $pa<x$ and $x=(p+1)a-c$ for some $c$ (possibly negative), then
$x\in\mathcal{A}(a)$ implies $c\leq p+1$.
\end{itemize}
\end{fact}

Now we can state the main theorem of this section:

\begin{theorem}
\label{theoremdiffaa}
If $A=(1,a_1,\ldots,a_k)$ is orderly then 
$$a_{j}-a_{i} \in \mathcal{A}(a_1)$$
for all $0\leq i<j\leq k$.
\end{theorem}

{\bf Proof.} If not then by property (3) above there exists a number $j$
for which $a_j-a_{j-1}\not\in\mathcal{A}(a_1)$, which is equivalent to
$$pa_1+1\leq a_j-a_{j-1}\leq (p+1)a_1-(p+2)$$
for some $p$. Among all pairs $(p,j)$ for which these inequalities hold
let us choose the lexicographically smallest one. Comparing the leftmost and
rightmost expressions in this double inequality yields $a_1\geq p+3$, hence
$a_1\geq 4$ and $1\leq p \leq a_1-3$.

We have $\opt_A(a_{j-1}+(p+1)a_1)\leq p+2$ and
$$a_j+(p+2)\leq a_{j-1}+(p+1)a_1\leq a_j+a_1-1$$
hence $\grd_{(1,\ldots,a_j)}(a_{j-1}+(p+1)a_1)\geq p+3$. It follows
that $a_{j+1}\leq a_{j-1}+(p+1)a_1$. Then
$$a_{j+1}-a_j\leq a_{j-1}+(p+1)a_1 -a_j \leq a_1-1.$$
By Proposition \ref{lemmadiffa1-1} all these inequalities must in fact be
equalities. In other words:
\begin{center}
\begin{tabular}{c}
$a_{j+1}=a_{j-1}+(p+1)a_1,$\\
$a_j=a_{j-1}+pa_1+1.$
\end{tabular}
\end{center}
Choose the largest $l$ for which $a_{l+1}-a_l\leq a_{j-1}+(p-1)a_1+2$ (such
$l$'s exist, for instance $a_{j+1}-a_{j}=a_1-1$ is sufficiently small). By
maximality of $l$ we have $a_{l+2}-a_{l+1}\geq a_{j-1}+(p-1)a_1+3$ (it
is possible that $a_{l+2}=\infty$). Observe that 
$$a_{l+2}-a_l = (a_{l+2}-a_{l+1})+(a_{l+1}-a_l)\geq a_{j-1}+(p-1)a_1+3 +
a_1-1=a_{j-1}+pa_1+2=a_j+1$$
and
$$a_{l+1}-a_l\leq a_{j-1}+(p-1)a_1+3 = a_j-a_1+2$$
which means that
$$a_{l+1}+a_1-2\leq a_l+a_j<a_{l+2}.$$
This eventually implies that $a_l+a_j=a_{l+1}+a_r$ for some $1\leq r<j$.
The rest of the proof depends on the possible locations of $a_l+a_{j-1}$.
 
If $a_l+a_{j-1}>a_{l+1}$ then the same argument yields an index $s<j-1$ for
which $a_l+a_{j-1}=a_{l+1}+a_s$. In that case
$a_r-a_s=a_j-a_{j-1}\not\in\mathcal{A}(a_1)$. By properties (3) and (2) of
$\mathcal{A}(a_1)$ there exist numbers $s<r'\leq r$ and $p'$ for which
$$p'a_1+1\leq a_{r'}-a_{r'-1}\leq (p'+1)a_1-(p'+2).$$
The inequality $a_{r'}-a_{r'-1}\leq a_j-a_{j-1}$ implies $p'\leq p$. The pair
$(p',r')$ is lexicographically smaller that $(p,j)$, which is a contradiction
since the latter was chosen to be minimal.

If, on the other hand, $a_l+a_{j-1}=a_{l+1}$,
then $a_l+a_j=a_l+a_{j-1}+pa_1+1=a_{l+1}+pa_1+1$, which means that $a_r=pa_1+1$.
Then $a_r-a_1=(p-1)a_1+1\not\in \mathcal{A}(a_1)$, contradicting the
minimality of $(p,j)$ be the same argument as above.

Therefore we are left with the case $a_l+a_{j-1}<a_{l+1}$. The number $a_r$
satisfies
$$a_r=a_l+a_j-a_{l+1}<a_l+a_j-(a_l+a_{j-1})=a_j-a_{j-1}=pa_1+1.$$
Since $a_r-a_1<(p-1)a_1+1$, the minimality of $p$ implies that
$a_r-a_1\in\mathcal{A}(a_1)$. It means that $a_r=qa_1-q'$ for some $q\leq p$ and
$0\leq q'<q$.

Next we are going to show that $a_{l+1}-(a_l+a_{j-1})\not\in\mathcal{A}(a_1)$.
Observe that
$$a_{l+1}-(a_l+a_{j-1})=(a_l+a_j-a_r)-a_l-a_{j-1}=pa_1+1-a_r=(p-q)a_1+(1+q').$$
which is more than $(p-q)a_1$, while at the same time it equals:
$$(p-q+1)a_1-(a_1-1-q')$$
with $a_1-1-q'>(p+2)-1-q'=p+1-q'>p-q+1$. By property (2) the number
$a_{l+1}-(a_l+a_{j-1})$ does not belong to $\mathcal{A}(a_1)$.

\begin{figure}[!h]
\epsfbox{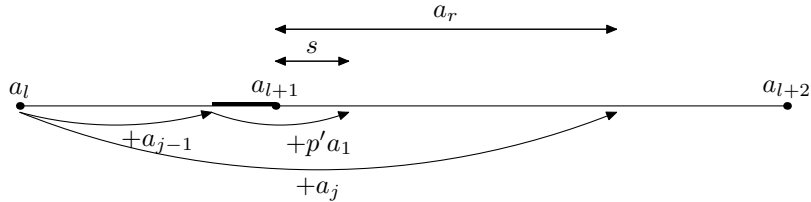}
\caption{The last case of the proof. The length of the bold interval is not
in $\mathcal{A}(a_1)$.}
\end{figure}

Now let us choose the least $p'$ for which $a_l+a_{j-1}+p'a_1\geq a_{l+1}$. In
this case
$$a_l+a_{j-1}+p'a_1<a_{l+1}+a_1\leq a_{l+1}+a_r=a_l+a_j.$$
Obviously $\opt_A(a_l+a_{j-1}+p'a_1)\leq p'+2$. On the other hand, the greedy
decomposition of $a_l+a_{j-1}+p'a_1$ is $a_{l+1}+s\cdot 1$, where
$s=a_l+a_{j-1}+p'a_1-a_{l+1}$. By optimality
$$s+1\leq p'+2$$
so $s\leq p'+1$. On the other hand, we have already proved that
$p'a_1-s=a_{l+1}-(a_l+a_{j-1})\not\in\mathcal{A}(a_1)$, so $s\geq p'+1$. Finally
we have $s=p'+1$.

To end the proof we compute $a_r-a_1$ in terms of $p,p'$ and $a_1$:
\begin{center}
\begin{tabular}{c}
$a_r-a_1=a_l+a_j-a_{l+1}-a_1=a_l+(a_{j-1}+pa_1+1)-(a_l+a_{j-1}
+p'a_1-s)-a_1=$\\
$=(p-p'-1)a_1+(s+1)=(p-p'-1)a_1+(p'+2)=$\\
$=(p-p')a_1-(a_1-p'-2)$
\end{tabular}
\end{center}
Since $a_r-a_1\in\mathcal{A}(a_1)$, by property (4) we obtain
\begin{center}
\begin{tabular}{c}
$a_1-p'-2\leq p-p',$\\
$a_1<p+3.$
\end{tabular}
\end{center}
This contradiction ends the proof.\qed

As an immediate corollary we obtain the theorem announced in the introduction:

\begin{theorem}
\label{theoremprefixorderly} 
 If $A=(1,a_1,\ldots,a_k)$ is orderly then for any $2\leq l\leq
k$ the currency $(1,a_1,a_l)$ is also orderly. In particular the currency
$(1,a_1,a_2)$ is orderly.
\end{theorem}

{\bf Proof}. If $A$ is orderly then by Theorem \ref{theoremdiffaa} we have
$a_l-a_1\in\mathcal{A}(a_1)$. By Proposition \ref{lemma3orderly} this is
sufficient for $(1,a_1,a_l)$ to be orderly.\qed

\section{Short currencies}
\label{section4}

Theorems \ref{theoremprefixorderly} and \ref{onept} allow us to give a
complete characterization of all orderly currencies with at most 5 coins. The
currencies with 1, 2 and 3 coins have already been discussed. Here we
concentrate on the cases of 4 and 5 coins. Following \cite{CCS} call a currency
$A=(1,a_1,\ldots,a_k)$ \emph{totally orderly}\footnote{Also called normal in
\cite{TienHu}.} if every prefix sub-currency of the form $(1,a_1,\ldots,a_l)$ is
orderly for $l=0,\ldots,k$.

\begin{proposition}
\label{lemma4}
The currency $A=(1,a_1,a_2,a_3)$ is orderly if and only if it is totally
orderly.
\end{proposition}

\begin{proposition}
\label{lemma5}
The currency $A=(1,a_1,a_2,a_3,a_4)$ is orderly if and only if 
\begin{itemize}
\item (1) either $(1,a_1,a_2,a_3,a_4)=(1,2,a,a+1,2a)$ for some $a\geq 4$, in
which case $(1,a_1,a_2,a_3)$ is not orderly,
\item (2) or $A$ is totally orderly.
\end{itemize}
\end{proposition}

{\bf Remark.} The conditions given in the above propositions are efficiently
computable, since it can be quickly checked if a currency is totally orderly (as
opposed to checking whether it is just orderly). One simply repeats the one-point
theorem test with for longer and longer prefixes; see also \cite{CCS}.

{\bf Proof of Propositions \ref{lemma4} and \ref{lemma5}.} The one-point
theorem, together with Theorem \ref{theoremprefixorderly} covers Proposition
\ref{lemma4} and case (2) of Proposition \ref{lemma5}. 

It remains to show that all orderly currencies $(1,a_1,a_2,a_3,a_4)$ in which
the sub-currency $(1,a_1,a_2,a_3)$ is disorderly are of the form (1) from
Proposition \ref{lemma5}. Let $m=\lceil a_3/a_2\rceil$.

The triple $(1,a_1,a_2)$ is orderly by Theorem \ref{theoremprefixorderly}. By
the one-point theorem $ma_2$ is a counterexample for $(1,a_1,a_2,a_3)$, hence
$a_4\leq ma_2$. Both values $a_3+a_3$ and $a_3+a_2$ exceed $ma_2$, so they
exceed $a_4$, so by optimality there must exist $i<j\leq 2$ for which:
\begin{center}
\begin{tabular}{c}
$a_3+a_2=a_4+a_i,$\\
$a_3+a_3=a_4+a_j.$
\end{tabular}
\end{center}
Subtracting these equations we get
$$a_3-a_2=a_j-a_i<a_j\leq a_2$$
which in turn gives $a_3<2a_2$. That means $m=2$. 

There are two cases to consider:

{\bf $j=2$.} Then $a_3-a_2=a_2-a_i$, so $2a_2=a_3+a_i$, which contradicts the fact that $(1,a_1,a_2,a_3)$ is disorderly.

{\bf $j=1$.} Then $i=0$ and previous equations take the form:
\begin{center}
\begin{tabular}{c}
$a_3+a_2=a_4+1,$\\
$a_3+a_3=a_4+a_1.$
\end{tabular}
\end{center}
The following computation
$$a_4+1=a_3+a_2>2a_2=ma_2\geq a_4$$
implies
$$a_4+1=a_3+a_2=2a_2+1.$$
Setting $a_2=a$ we get $a_3=a+1$, $a_4=2a$ and
$a_1=2a_3-a_4=2$.

The routine check that $(1,2,a,a+1,2a)$ is orderly resembles the technique used
in the proof of case 2 of Theorem \ref{onept} and is left to the reader.
For $a\geq 4$ the sub-currency $(1,2,a,a+1)$ is disorderly.\qed

Attempts to continue similar reasoning with longer coinage systems encounter a
serious problem, because the applicability of the one-point theorem is limited.
More precisely, the ``intermediate'' currencies may not be orderly even if $A$
is orderly as we see from part (1) of Proposition \ref{lemma5}. We shall return
to these matters in section \ref{section6}.

\section{Investigating differences, part II}
\label{section5}

In the previous sections we were discussing relation of the distances
$a_j-a_i$ and the value of $a_1$. Here we shall extend some of this to further
coins. Note that Proposition \ref{lemmadiffa1-1} may be interpreted as follows:
if some difference $a_j-a_i$ belongs to the interval $(1,a_1)$, then it must be
necessarily equal $a_1-1$. We are interested in the possible values of
$a_j-a_i$ in the cases when this difference belongs to $(a_{m-1},a_m)$. {\bf
Throughout this section we always assume that $A=(1,a_1,\ldots,a_k)$ 
is orderly.} The key results of this sections are Corollary \ref{possibletoa2}
and Theorem \ref{theorembigdiff}.

\begin{lemma}
\label{aux1}
If 
$$a_m-a_{l+1}<a_j-a_i<a_m-a_l$$
for some $i<j$, $l<m$, then
$$a_{j+1}\leq a_i+a_m.$$
\end{lemma}

{\bf Proof.} We have $a_j+a_l<a_i+a_m<a_j+a_{l+1}$. If there was no new coin
between $a_j$ and $a_i+a_m$ then there would be no greedy decomposition of
$a_i+a_m$ in two steps.\qed

\begin{lemma}
\label{aux2}
There are no numbers $0\leq i<j\leq k$ and $1\leq m\leq k$ that satisfy 
$a_{m-1}\leq a_j-a_i < a_m-a_{m-1}$.
\end{lemma}

{\bf Proof.} Suppose the contrary and let $(j,i,m)$ be some triple
satisfying the above inequalities, such that $j$ is the least possible. If
$i=0$ then $a_{m-1}\leq a_j-1<a_m-a_{m-1}<a_m$, hence $j=m$, which in turn
implies $a_{m-1}<1$, but this is not possible. 

Therefore $i\geq 1$ and we are free to choose the largest index $l$ for which
$a_l-a_{l-1}<a_i$. If $a_{l-1}+a_j<a_{l+1}$ then by the standard argument we
obtain a contradiction with the minimality of $j$. Hence $a_{l-1}+a_j\geq
a_{l+1}$. It follows that
$$a_l-a_{l-1}=(a_{l+1}-a_{l-1})-(a_{l+1}-a_l)\leq a_j-a_i<a_m-a_{m-1}.$$
By Lemma \ref{aux1} it follows that $a_{l+1}-a_{l-1}\leq a_m$, so
$a_i<a_{l+1}-a_{l-1}\leq a_m$. In effect $i\leq m-1$. At the same time we
also have $a_j>a_{m-1}$, so $j\geq m$. All this implies
$$a_j-a_i\geq a_m-a_{m-1}.$$
This contradiction ends the proof.\qed

\begin{lemma}
\label{aux3}
Let $m\geq 2$. If the difference $a_j-a_i$ belongs to the interval
$[a_m-a_1,a_m-1]$ then it can only be one of the numbers $a_m-a_1$, $a_m-a_1+1$
and $a_m-1$. 
\end{lemma}

{\bf Proof.} Suppose that
$$a_m-a_1<a_j-a_i<a_m-1.$$
From Lemma \ref{aux1} we get $a_{j+1}\leq a_i+a_m<a_j+a_1$. In this case
Proposition \ref{lemmadiffa1-1} implies $a_{j+1}=a_j+a_1-1$. 
Moreover, we have
$$a_j+2\leq a_i+a_m\leq a_j+a_1-1=a_{j+1}.$$
Hence $a_i+a_m=a_{j+1}$ (otherwise the amount $a_i+a_m$ would not have a greedy
decomposition in two steps). Eventually we get
$$a_j-a_i=(a_{j+1}-a_1+1)-(a_{j+1}-a_m)=a_m-a_1+1.$$\qed

\begin{lemma}
\label{aux4}
If $a_1<a_2-a_1+1<a_2-1$ then the value $a_2-a_1+1$ cannot be
attained by any of the differences $a_j-a_i$. 
\end{lemma}

{\bf Proof.} First note that the given inequalities imply $a_1\geq 3$.
Suppose that $j$ is the minimal number for which there exists an $i$ such that
$a_j-a_i=a_2-a_1+1$. Clearly $i\geq 2$. From the proof of Lemma \ref{aux3} we
know that
$$a_{j+1}=a_i+a_2=a_j+a_1-1.$$
Let $m$ be the maximal index for which $a_m-a_{m-1}<a_i$. Then $a_{m+1}-a_m\geq a_i$. 

If $a_{m-1}+a_j<a_{m+1}$ then considering the amounts $a_{m-1}+a_i$ and $a_{m-1}+a_j$ and their
greedy decompositions we obtain a contradiction with the minimality of $j$ in the usual way. Hence
we may assume that
$$a_{m-1}+a_j\geq a_{m+1}.$$

If $a_{m-1}+a_j=a_{m+1}$ then consider the amount $a_{m-1}+a_{j+1}$. It
satisfies
$$a_{m-1}+a_{j+1}=a_{m+1}+a_1-1<a_{m+1}+a_2\leq a_{m+1}+a_i\leq a_{m+2}.$$
Since $a_1-1\geq 2$ this amount cannot be greedily decomposed in two steps, so
we have a contradiction, which means that
$$a_{m-1}+a_j\geq a_{m+1}+1$$
which in turn implies
$$a_m-a_{m-1}=(a_{m+1}-a_{m-1})-(a_{m+1}-a_m)\leq a_j-1-a_i=a_2-a_1.$$
We know from the previous lemmas that in this case the only possible values of the difference
$a_m-a_{m-1}$ are $a_2-a_1$, $a_1$ and $a_1-1$. Let us investigate these cases separately.

{\bf Case 1.} $a_m-a_{m-1}=a_2-a_1$. Then $a_{m+1}\geq a_m+a_i$ and
$$a_{m+1}\leq a_{m-1}+a_j-1=a_m-a_2+a_1+a_j-1=a_m+a_i$$
hence $a_{m+1}=a_m+a_i=a_{m-1}+a_j-1$. Now consider the amount $a_m+a_j$. It
satisfies
$$a_m+a_j=a_m+a_i+a_2-a_1+1=a_{m+1}+a_2-a_1+1<a_{m+1}+a_2\leq a_{m+1}+a_i\leq a_{m+2}$$
so it could be decomposed greedily in two steps only if $a_2-a_1+1$ was one of
the coins, which is not true by the assumptions of the lemma.

{\bf Case 2.} $a_m-a_{m-1}=a_1$. Now consider the amount $a_{m-1}+a_2$:
$$a_m<a_{m-1}+a_2=a_m+(a_2-a_1)<a_m+a_i\leq a_{m+1}.$$
This amount can only be decomposed optimally if $a_2-a_1$ is a coin. Since
$a_1\geq 3$, by Proposition \ref{lemmanot1} we have $a_2-a_1\neq 1$. Therefore
$a_2-a_1=a_1$ and we have
$$a_m-a_{m-1}=a_1=a_2-a_1$$
and the argument from case 1 can be repeated.

{\bf Case 3.} $a_m-a_{m-1}=a_1-1$. An exact repetition of case 2 shows that in
this case $a_2-a_1+1$ would have to be one of the coins. However, this
possibility is excluded by the assumptions of our lemma. \qed

The results from this section, together with Proposition \ref{lemmadiffa1-1} can
be used to characterize the set of possible values of $a_j-a_i$ which fit in the
interval $(1,a_2)$. For a currency $A$ let $S(A)=\{a_j-a_i: 0\leq
i<j\leq k\}$.

\begin{corollary} 
\label{possibletoa2}
For an orderly currency $A=(1,a_1,a_2,\ldots,a_k)$
\begin{itemize}
\item[(a)] we always have 
$$S(A)\cap (1,a_1)\subset\{a_1-1\}$$
$$S(A)\cap (a_1,a_2)\subset\{a_2-a_1,a_2-1\}$$
\item[(b)] if $a_2=2a_1-1$ or $a_2=2a_1$ then $S(A)\cap
(1,a_2)\subset\{a_1-1,a_1,a_2-1\}$
\item[(c)] if $a_2>2a_1$ then $S(A)\cap (1,a_2)=\{a_1-1,a_2-a_1,a_2-1\}$
\end{itemize}
\end{corollary}

{\bf Proof.} Property (a) is just a restatement of Proposition
\ref{lemmadiffa1-1} and Lemmas \ref{aux2}, \ref{aux3} and \ref{aux4}.

By Theorem \ref{theoremdiffaa} there are no other possible values of $a_2$
except of those in (b) and (c). In both cases, if $a_j-a_i < a_1$, then
Proposition \ref{lemmadiffa1-1} applies. 

In case (c) $a_1<a_2-a_1<a_2-a_1+1\leq a_2-1$ and an application of Lemmas
\ref{aux3} and \ref{aux4} proves that our theorem enumerates all possible
elements of $S(A)\cap (1,a_2)$. Of course all the given values are attained, so
in (c) we are free to use equality rather than inclusion.

In case (b) the difference $a_1$ may or may not be attained (consult the
currencies $(1,3,5)$ and $(1,3,5,8,10,15)$). Once again one needs to combine
the before-mentioned lemmas; we omit the details.\qed

{\bf Remark.} Corollary \ref{possibletoa2} and Theorem \ref{theoremdiffaa} give
two independent conditions that must be satisfied by orderly currencies. For
instance every three-coin currency satisfies Corollary \ref{possibletoa2}, but
not necessarily Theorem \ref{theoremdiffaa} (it is also easy to imagine more
complicated examples of this kind). On the other hand, the currency $(1,3,7,12)$
satisfies Theorem \ref{theoremdiffaa}, but $12-7=5\not\in\{7-3,7-1\}$, so part
(a) of Corollary \ref{possibletoa2} is violated.

Our last theorem in this section will be important in section
\ref{section7}. It can roughly be stated as ``if some two consecutive
differences are large, then the subsequent differences must also be large''.

\begin{theorem}
\label{theorembigdiff}
Suppose $(1,a_1,\ldots,a_k)$ is orderly, $m\geq 2$ and
$$a_{m-1}>2a_{m-2}, \ a_{m}>2a_{m-1}.$$
Then for every $t\geq m$ we have $a_{t+1}-a_t\geq a_m-a_{m-1}$.
\end{theorem}

{\bf Proof.} Suppose, on the contrary, that $a_{t+1}-a_t< a_m-a_{m-1}$ for some
$t\geq m$, and let $t$ be the smallest index with these properties. Choose $s$
as the largest index for which
$$a_{s+1}-a_s<a_t-a_{m-2}$$
(such numbers $s$ exist; for instance $s=t-1$ satisfies this inequality).
Note that by maximality of $s$ we have $a_{s+2}-a_{s+1}\geq a_t-a_{m-2}$
(possibly $a_{s+2}=\infty$) and $a_{s+3}-a_{s+2}\geq a_t-a_{m-2}$ (if
$a_{s+2}<\infty$). The proof is split into two cases.

{\bf Case 1.} $a_s+a_{t+1}<a_{s+2}$. With this assumption we have
$$a_{s+1}<a_s+a_t<a_s+a_{t+1}<a_{s+2}$$
so there exist indices $r,l$ such that
\begin{center}
\begin{tabular}{c}
$a_s+a_t=a_{s+1}+a_r,$\\
$a_s+a_{t+1}=a_{s+1}+a_l,$
\end{tabular}
\end{center}
with $r<l\leq t$. This implies 
$$a_l-a_{l-1}\leq a_l-a_r=a_{t+1}-a_t<a_m-a_{m-1}.$$
Since $l-1< t$ and $t$ was chosen to be minimal with respect to the condition
$t\geq m$ and the above inequality, we obtain $l-1<m$. Since $l=m$ does not
satisfy the above inequality, we have $l\leq m-1$ and $r\leq m-2$, but then
$$a_{s+1}-a_s=a_t-a_r\geq a_t-a_{m-2}$$
contradicting the choice of $s$. This completes the first case of the proof.

{\bf Case 2.} Now suppose $a_s+a_{t+1}\geq a_{s+2}$. We are going to prove the
following sequence of inequalities:
\begin{center}
\begin{tabular}{lc}
(1) & $a_{s+1}-a_s>a_{m-2}$\\
(2) & $a_{s+2}-a_s>a_m$\\
(3) & $a_{s+1}-a_s<a_m$\\
(4) & $a_{s+1}-a_s\geq a_m-a_{m-1}$\\
(5) & $a_{s+2}<a_{s+1}+a_t<a_{s+1}+a_{t+1}<a_{s+3}$
\end{tabular}
\end{center}

(1): We always have
$$a_{s+2}-a_s>a_{s+2}-a_{s+1}\geq a_t-a_{m-2}\geq a_m-a_{m-2}>a_{m-1}.$$
If we also had $a_{s+1}-a_s\leq a_{m-2}$ then
$$a_{s+1}\leq a_s+a_{m-2}<a_s+a_{m-1}<a_{s+2}.$$
As usually, it means that $a_{m-1}-a_{m-2}=a_l-a_r$ for some $r<l\leq m-2$ or
$a_{m-1}-a_{m-2}=a_l$ for $l\leq m-2$. In either case $a_{m-1}-a_{m-2}\leq
a_{m-2}$, contradicting the assumptions of the theorem. Therefore
$a_{s+1}-a_s>a_{m-2}$.

(2): This follows straight from (1) and the maximality of $s$:
$$a_{s+2}-a_s=(a_{s+2}-a_{s+1})+(a_{s+1}-a_s)>a_t-a_{m-2}+a_{m-2}=a_t\geq a_m.$$

(3): Since we assumed $a_{s+2}-a_s\leq a_{t+1}$ for this case, we obtain, using
the properties of $s$ and $t$, that
$$a_{s+1}-a_s=(a_{s+2}-a_s)-(a_{s+2}-a_{s+1})\leq
a_{t+1}-(a_t-a_{m-2})<a_m-a_{m-1}+a_{m-2}<a_m.$$

(4): By (2) and (3) we have $a_{s+1}<a_s+a_m<a_{s+2}$, therefore
$a_s+a_m=a_{s+1}+a_r$ for some $r\leq m-1$. Finally
$$a_{s+1}-a_s=a_m-a_r\geq a_m-a_{m-1}.$$

(5): First note that by $a_t\geq a_m$ and $2a_{m-2}<a_{m-1}$ we obtain
$$a_{s+3}-a_{s+1}\geq 2(a_t-a_{m-2})=a_t+a_t-2a_{m-2}>
a_t+a_m-a_{m-1}>a_{t+1}.$$
Moreover, by (4) and the assumption $a_{s+2}-a_s\leq a_{t+1}$ we get
$$a_{s+2}-a_{s+1}=(a_{s+2}-a_s)-(a_{s+1}-a_s)\leq a_{t+1}-(a_m-a_{m-1})<a_t.$$
This ends the proof of (1)--(5).
\begin{figure}
\epsfbox{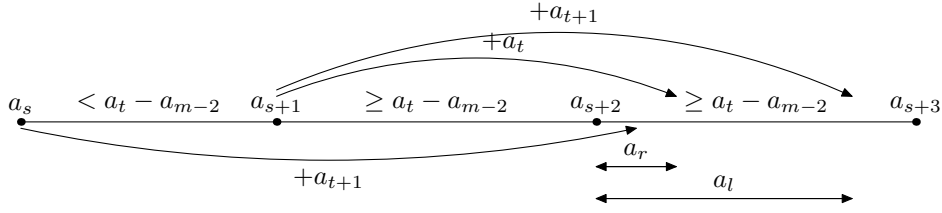}
\caption{The situation in case 2. in Theorem \ref{theorembigdiff}.}
\end{figure}
Now (5) implies the existence of $r<l\leq t$ such that
\begin{center}
\begin{tabular}{c}
$a_{s+1}+a_t=a_{s+2}+a_r,$\\
$a_{s+1}+a_{t+1}=a_{s+2}+a_l.$
\end{tabular}
\end{center}
As a consequence of these formulae we obtain the inequality
$$a_r=a_t-(a_{s+2}-a_{s+1})\leq a_t-(a_t-a_{m-2})=a_{m-2}, \ \textrm{hence }
r\leq m-2,$$
which in turn implies
$$a_l=(a_{t+1}-a_t)+a_r<a_m-a_{m-1}+a_{m-2}<a_m, \ \textrm{hence } l\leq m-1.$$
Combining this, we get
$$a_{s+1}-a_s=a_{s+2}-a_s+a_l-a_{t+1}\leq a_{t+1}+a_l-a_{t+1}=a_l\leq a_{m-1}.$$
However, by (4) $a_{s+1}-a_s\geq a_m-a_{m-1}>a_{m-1}$, so we have a
contradiction which ends the proof of case 2, and the whole theorem.\qed

\section{$+/-$-classes}
\label{section6}

If $A=(1,a_1,\ldots,a_k)$ is orderly then some prefix sub-currency, i.e. a
currency of the form $A'=(1,a_1,\ldots,a_l)$ with $l<k$ might not be orderly
(for instance, $(1,2,a,a+1,2a)$ is orderly, but $(1,2,a,a+1)$ is not for $a\geq
4$, as in Proposition \ref{lemma5}).
This situation was still quite manageable in the case of 5 coins, but it gets
more and more complicated as the number of coin increases, thus making inductive
analysis (possibly using the one-point theorem) impossible. 

To describe the prefix currencies we introduce the notion of {\it
$+/-$-classes}. To every currency $A=(1,a_1,\ldots,a_k)$ we may assign a
pattern of $k+1$ signs {\tt +} and {\tt -}, defined as follows: the $l$-th
symbol of the pattern ($l=0,\ldots,k$) is {\tt +} if the prefix currency
$(1,a_1,\ldots,a_l)$ is orderly and {\tt -} in the opposite case. A {\it
$+/-$-class} is the set of all currencies corresponding to a given
$+/-$-pattern. For instance, the pattern \verb?++++?$\ldots$\verb?+++?
corresponds to totally orderly currencies. Another well-described example
is the $+/-$-class given by the pattern \verb?+++-+? --- it consists precisely of the
currencies $(1,2,a,a+1,2a)$ with $a\geq 4$ (this is the consequence of 
Proposition \ref{lemma5}, since an orderly $5$-coin currency which is not
totally orderly satisfies part (1) of that proposition).

The $+/-$-patterns that correspond to non-empty classes cannot be completely
arbitrary, for instance, if a pattern ends with a {\tt +} then it must begin
with {\tt +++} -- this is a consequence of Theorem \ref{theoremprefixorderly}.
The patterns beginning with {\tt +++} and ending with {\tt +} will be called
{\it proper}. Mysteriously, some proper patterns describe empty classes. Here is
a sample proposition of this sort:

\begin{proposition}
\label{emptyclass}
The $+/-$-class described by the pattern \verb?+++-+-+? is empty.
\end{proposition}

{\bf Proof.} Suppose that $A=(1,a_1,a_2,a_3,a_4,a_5,a_6)$ is a coinage system in
the class \verb?+++-+-+?. By case (1) of Proposition \ref{lemma5} we know that
in fact $A$ is of the form 
$$(1,2,a,a+1,2a,a_5,a_6)$$
for some $4\leq a$, $2a< a_5< a_6$. By the one-point theorem some multiple of
$2a$ is a counterexample for $(1,2,a,a+1,2a,a_5)$. Extending this by $a_6$ must
fix this problem, hence
$$a_6-a_5<2a.$$
Since $A$ is orderly, there exist numbers $r, s$ such that:
\begin{center}
\begin{tabular}{c}
$a_5+a_5=a_6+a_r,$\\
$a_5+2a=a_6+a_s,$
\end{tabular}
\end{center}
with $a_r\leq 2a$, $a_s\leq a+1$, $1 \leq a_s<a_r$. Subtracting the two
equations yields $a_5-2a=a_r-a_s$. Possible differences $a_r-a_s$ ($0\leq
s<r\leq 4$) form the set
$$\{1,a-2,a-1,a,2a-2,2a-1,2a\}$$
so the possible values of $a_5$ are $2a+1,3a-2,3a-1,3a,4a-2,4a-1,4a$.

The values $3a-1,3a,4a-2,4a-1,4a$ can be excluded from this set, since then
$(1,2,a,a+1,2a,a_5)$ would be orderly, which can be checked easily by the
one-point theorem (the ``suspected'' amount to be tested for optimality is
$4a$).

Therefore we are left with $a_5\in\{2a+1,3a-2\}$.

If $a_5=2a+1$ then the greedy algorithm for $(1,\ldots,a_5)$ fails to be optimal
already for $3a=2a+a$, hence $a_6\leq 3a$. On the other hand, all three numbers
$2a+2a$, $2a+(2a+1)$ and $(2a+1)+(2a+1)$ can be obtained with two coins, hence
$4a-a_6$, $4a+1-a_6$ and $4a+2-a_6$ must be three consecutive integers which are
coins, all less than $a_6$. This is only possible if $a_6=4a$, contradiction.

Now suppose that $a_5=3a-2$. Then for the number $2a+(a+1)=3a+1$ not to be a
counterexample we must have $3a-1\leq a_6\leq 3a+1$. If $a_6=3a-1$ then
$4a-2=(3a-2)+a=(3a-1)+(a-1)$ is a counterexample ($a-1$ is not a coin). If
$a_6=3a$ then the counterexample is $4a-1=(3a-2)+(a+1)=3a+(a-1)$ (reason as
before). Finally, if $a_6=3a+1$ then $4a=2a+2a=(3a+1)+(a-1)$ is the
counterexample.\qed

Of course, given a currency, we may recover its $+/-$-class in $O(k^4)$ time
simply by repeating Pearson's algorithm \cite{Pear} for each prefix
sub-currency. The reverse problem, to determine whether a given proper
$+/-$-pattern describes a non-empty $+/-$-class is actually much harder and we
have not been able to find any algorithm solving it.

From this point of view the most ``messy'' orderly currencies are those which
belong to the class determined by \verb?+++----?$\ldots$\verb?--+?. These
classes are indeed non-empty for $k=0,2\pmod{3}$. Their representatives for
$k=3l$ and $k=3l-1$, respectively, are
\begin{center}
\begin{tabular}{c}
$(1,2,\ 4,5,\ 7,8,\ldots,3l-2,3l-1,\ 3l+1,\ 3l+4,\ldots,6l-2),$\\
$(1,2,\ 4,5,\ 7,8,\ldots,3l-2,3l-1,\ 3l+2,\ 3l+5,\ldots,6l-4).$
\end{tabular}
\end{center}
On the other hand, there seem to be no coinage
systems of type \verb?+++----?$\ldots$\verb?--+? for $k=1\pmod{3}$, but we have
not been able to prove this.

\section{Classification of orderly sub-currencies}
\label{section7}

Every set $P=\{i_0,i_1,\ldots,i_l\}\subset\{0,1,\ldots,k\}$, where
$0=i_0<i_1<\ldots<i_l$ determines a sub-currency
$(a_{i_0},a_{i_1},\ldots,a_{i_l})$ of any currency $A=(1,a_1,\ldots,a_k)$. 
From Theorem \ref{theoremprefixorderly} we know that if $A$ is orderly then the
sub-currency determined by $P=\{0,1,l\}$ ($2\leq l\leq k$) is also orderly. Is
this just a lonely phenomenon, or could a similar theorem be proved for some
other sets $P$?

\begin{definition}
\label{hered}
The set $P$ of the form given above will be called \emph{hereditary} if the
following is true:
\begin{center}
for every orderly currency $A=(1,a_1,\ldots,a_k)$\\ the sub-currency
determined by $P$ is also orderly
\end{center}
\end{definition}

Let us enumerate some interesting classes of subsets of $\{0,1,\ldots,k\}$:

\begin{itemize}
\item[type 1:] the singleton set $\{0\}$
\item[type 2:] the sets $\{0,l\}$ for $1\leq l\leq k$
\item[type 3:] the sets $\{0,1,l\}$ for $2\leq l\leq k$
\item[type 4:] the sets $\{0,1,2,l\}$ for $4\leq l\leq k$
\item[type 5:] the full set $\{0,1,\ldots,k\}$
\end{itemize}

Note that $\{0,1,2,3\}$ is a peculiar exception: it is \emph{not} of type 4 (an
immediate example is $(1,2,a,a+1,2a)$ for $a\geq 4$ and its non-orderly
sub-currency $(1,2,a,a+1)$ determined by $\{0,1,2,3\}$).
 
We already know that sets $P$ of type 1, 2, 3 or 5 are hereditary. In this
section we shall prove that sets that are not specified in types 1--5 are not
hereditary.\footnote{To be precise, every set $P$ should always be thought of as
a subset of $\{0,1,\ldots,k\}$ for a certain $k$. In most cases $k$ will be
implicit, but to improve clarity we shall sometimes stress this connection by
writing $P\subset\{0,1,\ldots,k\}$.}
We also conjecture that all sets $P$ of type 4 are hereditary, and we prove this
conjecture under some mild additional assumptions. The general case remains
open.

Before proceeding with the elimination of non-hereditary subsets $P$ let us
make a few observations.

\begin{lemma}
For any $l\geq 3$ let $B_l$ denote the currency
$$B_l=(1,2,3,\ldots,l-1,2l-2,2l-1,4l-4)$$
where $a_l=2l-1$. Then $B_l$ is orderly of type \verb?+++?$\ldots$\verb?+-+?.
\end{lemma}

{\bf Proof.} The prefix currency $(1,2,3\ldots,l-1)$ is clearly of type
\verb?+++?$\ldots$\verb?++?. Extending this by $2(l-1)$ we get an orderly
currency by the one-point theorem. The next prefix, ending in $2l-1$ is not
orderly since $2\cdot 2(l-1)=4l-4$ is the smallest counterexample. The complete
currency is orderly which can be proved easily by the techniques from the proof
of Theorem \ref{onept}.\qed

\begin{lemma}
For any $m>l\geq 2$ and $p\geq 1$ let $A_{l,m}(p)$ denote the currency:
$$A_{l,m}(p) = (a_0,a_1,a_2,\ldots,a_{l-1},a_l,a_{l+1},a_{l+2},\ldots,a_{m})=$$
$$(1,2,3,\ldots,l,pl,(2p-1)l,(3p-2)l,\ldots,((m-l+1)p-(m-l))l)$$
where $a_l=pl$. This currency is orderly. Moreover, if $p>m-l$ then $\lceil
a_m/a_l \rceil = m-l+1$.
\end{lemma}

{\bf Proof.} The given currency is in fact of type
\verb?++++?$\ldots$\verb?+++?, which can be verified inductively by the
one-point theorem: to check that $(1,a_1,\ldots,a_{l+i})$ is orderly for $i \geq
1$ it suffices to observe that 
$$2a_{l+i-1} = 2(pi-(i-1))l= (p(i+1)-i)l + (p(i-1)-(i-2))l=a_{l+i}+a_{l+i-2}.$$
To prove the last statement note that
$$a_{m}=((m-l+1)p-(m-l))l<(m-l+1)pl=(m-l+1)a_l$$
and, if $p>m-l$:
$$a_{m}=((m-l+1)p-(m-l))l = (m-l)pl + l(p-(m-l)) > (m-l)a_l.$$\qed

\begin{lemma}
\label{obs3}
An orderly currency may be extended by any multiple of its
highest coin and the resulting currency will be orderly.
\end{lemma}

{\bf Proof.} A trivial consequence of the one-point theorem.\qed

The last observation will be used in the following way: suppose we want to
prove that some set $P=\{i_0,i_1,\ldots,i_l\}\subset\{0,1,\ldots,k\}$ is not
hereditary. First we find a shorter orderly currency $A'=(1,a_1,\ldots,a_r)$,
such that the sub-currency determined by
$P'=\{i_0,i_1,\ldots,i_{r'}\}\subset\{0,1,\ldots,r\}$ is not orderly (here
$r'< r\leq k$) and $i_{r'+1}>r$ or $r'=l$. Let $c$ be any counterexample for
this sub-currency and let $m$ be any number for which $ma_{r}>c$. Then the
currency $$A=(1,a_1,\ldots,a_r,ma_r,2ma_r,\ldots,(k-r)ma_r)$$
is orderly (Lemma \ref{obs3}) and its sub-currency determined by $P$ is not,
since all the added coins are too large to fix the problem with $c$ (the exact
form of $P\setminus P'$ is actually immaterial, it is important that its
smallest element is at least $r+1$).

\begin{theorem}
The sets $P$ not of the form 1, 2, 3, 4 or 5 are not hereditary.
\end{theorem}

{\bf Proof.} Let $P=\{i_0,\ldots,i_s\}\subset\{0,\ldots,k\}$, $i_0=0$, be such a
set. Let $r$ be the largest index for which $i_r=r$ (i.e. $\{0,\ldots,r\}\subset
P$, $r+1\not\in P$). We shall consider a few cases:

{\bf Case $3\leq r<k$}. Here we employ the orderly currency $B_r$. Its
sub-currency $(1,a_1,\ldots,a_r)$ is not orderly. If $r=k-1$ then we are done,
while for $r<k-1$ we must expand $B_r$ to an orderly currency with $k+1$
coins in the standard way described earlier. The resulting currency will
have a disorderly sub-currency determined by $P$.\

{\bf Case $r=2$}. In this case $|P|\geq 5$, since otherwise $P$ would be of
the form $\{0,1,2\}$ or $\{0,1,2,l\}$ for some $l\geq 4$ and these sets are of
type 3 and 4, respectively. Denote $l=i_3\geq 4$, $m=i_4$ and consider the
currency $A_{l,m}(p)$ with $p>m-l$. Its sub-currency
$$(1,2,3,a_l,a_m)$$
is not orderly since the amount
$$\lceil a_m/a_l\rceil a_l = (m-l+1)a_l$$
paid greedily splits into the coin $a_m$ and some of the coins $1,2,3$, thus
requiring at least
$$1+\frac{(m-l)l}{3} > 1+(m-l)$$
coins, which is more than if it was paid with $m-l+1$ copies of $a_l$. Now it
suffices to expand this currency to a currency with $k+1$ coins as previously.

{\bf Case $r=1$}. Then $|P|\geq 4$, since otherwise $P$ would be of the form
$\{0,1,l\}$, which is of type 3. Let $l=i_2\geq 3$ and $m=i_3$ and consider the
currency $A_{l,m}(p)$ with $p>m-l$. The sub-currency $(1,2,a_l,a_m)$
is not orderly for the same reason as previously: the amount 
$\lceil a_m/a_l\rceil a_l = (m-l+1)a_l$
must be paid greedily with at least
$1+\frac{(m-l)l}{2} > 1+(m-l)$
coins and the proof follows.

{\bf Case $r=0$}. Clearly $|P|\geq 3$, since sets of the form $\{0\}$ and
$\{0,l\}$ are of type 1 and 2. Let $l=i_1\geq 2$ and $m=i_2$. Repeat the
same arguments with the currency $A_{l,m}(p)$ ($p>m-l$) and its sub-currency
$(1,a_l,a_m)$: this time the amount 
$\lceil a_m/a_l\rceil a_l = (m-l+1)a_l$
must be paid greedily with at least
$1+\frac{(m-l)l}{1} > 1+(m-l)$
coins.\qed

Sets $P$ of type 4 are the most peculiar ones. We believe they are also
hereditary; that is, we have the following:

\begin{conjecture}
\label{contype4}
If $A=(1,a_1,\ldots,a_k)$ is orderly, then the currency $(1,a_1,a_2,a_l)$ is
also orderly for every $4\leq l\leq k$.
\end{conjecture}

While this is not known to be true in general, we can prove this conjecture
under some mild additional conditions.

\begin{theorem}
\label{theorempartialconjecture}
Conjecture \ref{contype4} is true if we additionally assume that $a_2>2a_1$ and
$a_3>2a_2$.
\end{theorem}

{\bf Proof.} We shall verify that $(1,a_1,a_2,a_l)$ is orderly by
Proposition \ref{lemma4}. Let $m=\lceil a_l/a_2\rceil$.
By Theorem \ref{theorembigdiff} for every $l\geq 3$ we have the first
of the following inequalities:
$$a_{l+1}-a_l\geq a_3-a_2>a_2>ma_2-a_l.$$
It means that $a_{l+1}>ma_2$, so there is no new coin between $a_l$ and
$ma_2$, and the greedy decomposition of $ma_2$ with respect to $A$ involves
only the coins $1,a_1,a_l$. This justifies the first equality in the following
comparison:
$$\grd_{(1,a_1,a_2,a_l)}(ma_2)=\grd_{A}(ma_2)=\opt_{A}(ma_2)\leq
\opt_{(1,a_1,a_2,a_l)}(ma_2)$$
and by Proposition \ref{lemma4} the proof is complete.\qed





\section{Closing remarks and open problems}
\label{section8}

Throughout this paper we have proposed some possible approaches to the problem
of describing orderly coinage systems and their interesting properties. Some of
these techniques have enabled us to prove the most important results of
this paper, namely the structural theorems, like Theorem
\ref{theoremdiffaa} and Corollary \ref{possibletoa2}, or to give concise
descriptions of small systems. There is still quite a lot of work to be done in
the following areas:

\begin{itemize}
\item {\bf sub-currencies}: prove Conjecture \ref{contype4}, thus completing
the classification of orderly sub-currencies.

\item {\bf prefix sub-currencies}: invent an algorithm to decide whether a
given $+/-$--pattern describes a non-empty class or devise some other
properties of such $+/-$--patterns. Another interesting conjecture, to which we
have not found a counterexample, is:
\begin{conjecture}
If a $+/-$--class is non-empty, then it has a representative
$A=(1,a_1,\ldots,a_k)$ with $a_1=2$.
\end{conjecture}

\item {\bf differences}: can Corollary \ref{possibletoa2} be generalized? In
other words, what can be said about the differences $a_j-a_i$ that belong to
$(a_{m-1},a_m)$ for some $m$? Is it true that in general
$$S(A)\cap (a_{m-1},a_m)\subset \{a_m-a_{m-1},a_m-a_{m-2},\ldots,a_m-1\},$$
where $S(A)=\{a_j-a_i: 0\leq i<j\leq k\}$ for an orderly currency $A$? We
already know this is true for $m=1,2$. The lemmas from section \ref{section6}
provide some partial results in the general case as well.

\item {\bf extending}: Theorem \ref{theoremdiffaa}, Corollary
\ref{possibletoa2} and Conjecture \ref{contype4} can be thought of as {\it
obstructions} against extending: if a currency does not satisfy one of these
conditions then it cannot be extended to an orderly currency by appending new
coins of high denominations (higher than all the existing coins). What are the
other invariants of this sort? Is there an algorithm that decides if a currency
can be extended to an orderly one? Problems related to obstructions and
extending can also be found in \cite{CCS}.
\end{itemize}

{\bf Acknowledgements.} We are indebted to the referee, whose
valuable suggestions improved both the presentation and some technical aspects
of our paper. We also thank Lenore Cowen for pointing us to \cite{CCS}.


\end{document}